\documentclass[reqno]{amsart}
\usepackage[notref,notcite]{showkeys}

\begin{document}
\title
{A note on the Hodge conjecture}

\author
{Jorma Jormakka}

\email{jorma.o.jormakka@gmail.com}

\keywords{Hodge theory; differential forms;
algebraic geometry}

\begin{abstract}
The paper presents a counterexample to the Hodge conjecture.

\end{abstract}

\maketitle
\numberwithin{equation}{section}
\newtheorem{theorem}{Theorem}[section]
\newtheorem{lemma}[theorem]{Lemma}
\allowdisplaybreaks

\section{Introduction}

The Hodge conjecture is one of the better known open problems in mathematics 
and was chosen as one of the Millennium Prize problems by the Clay 
Mathematics Institute [1]. The formulation of the conjecture in [1] is: On
a projective non-sigular algebraic variety $C$, any Hodge class is a 
rational linear combination of classes $cl(\mathbb{Z})$ of algebraic
cycles. This seems to say that every Hodge class on an algebraic projective
complex manifold $M$ is a linear combination with rational coefficients of the 
cohomology classes obtained by the Poincar\'e duality from homology classes
of complex algebraic subvarieties of $M$.

The paper proposes a counterexample to the Hodge conjecure in the
the algebraic projective complex manifold
\begin{equation}
\begin{gathered}
M=\{(s_1:s_2:s_3:s_4:s_5:s_6:s_7:s_8)|s_1^4+s_2^4+s_3^4+s_4^4=0,\\
s_5^4+s_6^4+s_7^4+s_8^4=0,s_4^4+s_5^4=0\}
\end{gathered}
\end{equation}
$M$ is an algebraic subvariety of codimension 3 in the complex projective 
space $\mathbb{P}^7$ and a closed complex manifold of dimension 4. 
The manifold $M$ 
is composed as a product of two Fermat quadratic surfaces tied with the
third homogeneous polynomial equation in order to embed $M$ into a projective
space. As the Fermat quadratic surface is a $K3$ space, the manifold $M$ 
inherits a (2,0)-form that is never zero from the first copy of $K3$
and a similar (0,2)-form from the second copy of $K3$. 
It is shown in the presented paper 
that the wedge product of these forms is a (2,2)-form that
cannot be represented as a $\mathbb{Q}$-linear combination of 
cohomology classes of deriving from algebraic subvarieties of $M$. 
This form, having only one term 
$f(P)dz_1\wedge dz_2\wedge d\bar{z}_3\wedge d\bar{z}_4$, is rational 
(it is a linear combination of this type of forms with rational coefficients)
and thus a Hodge class. 

The study of what cohomology classes can be obtained by the Poincar\'e duality
from homology classes of submanifolds was started by Ren\'e Thom, but in real
manifolds. Thom's original paper is in French and rather difficult to read,
but we can look at [2] that is available on-line and find there the 
following statement (page 1 in [2]): All homology classes with integral 
coefficients of compact orientable differentiable manifolds of dimension 
$<10$ are realizable by submanifolds. As a 4-dimensional complex manifold can
be understood as an 8-dimensional real manifold, this statement, combined 
with the Poincar\'e duality, implies that every singular cohomology class 
of the
8-dimensional manifold can be realized as a linear combination with integer
coefficients of classes of real submanifolds. However, these real submanifolds
need not be complex submanifolds because a complex submanifold inherits 
the complex structure from the mother manifold. 

If $Z$ is a complex submanifold of $M$ of codimension $k$ then the 
inclusion map $i:Z\to M$
induces a cohomology class $[Z]$ that is in $H^{k,k}(M)$. In local coordinates
the cohomology class $[Z]$ is then a form of the type
\begin{equation}
\psi=f(z,\bar{z})dz^I\wedge d\bar{z}^J
\end{equation}
where $(z,\bar{z})=(z^1,\dots,z^n,\bar{z}^1,\dots,\bar{z}^n)$ are the local
coordinates. As $\psi\in H^{k,k}(M)$ there are the same number $k$ 
of indices in the multi-indices $I$ and $J$. This does not seem to say
very much. However, 
if we look more carefully 
we see that the form (1.2) has the property: if $z$ are the local 
coordinated of the submanifold,
the differentials $dz^j$ and $d\bar{z}^j$ always appear in pairs, as 
$dz^j\wedge d\bar{z}^j$. The reason is the following.

The tangent
space of a 2-dimensional complex manifold in a chosen base point $P$
is $\mathbb{C}^2$. There are two complex coordinates $z_1=x_1+iy_1$ 
and $z_2=x_2+iy_2$. (Let us use lower indices for a while 
as it is more natural.)
A 2-dimensional real plane in the tangent space of 
the chosen basepoint can be spanned by any linear combination of four
independent vectors $e_{x,1},e_{y,1},e_{x,2},e_{y,2}$ but a complex line
cannot separate $x_j$ and $y_j$. They are not vectors, they are numbers in
the field $\mathbb{C}$. Thus, a complex line through the basepoint 
can only be a linear combination $z=az_1+bz_2=c$, where $a,b$ and $c$ 
are complex 
numbers. Such a complex line defines a real plane as the  $x$ and $iy$ 
coordinates of $z=x+iy$ correspond to the vectors $e_x$ and $e_y$. 
In the local coordinates of the (real) mother manifold this (real) 
plane has the coordinates $e_x=ae_{x,1}+be_{x,2}$ and
$e_y=ae_{y,1}+be_{y,2}$. It can be expressed by using a wedge product
$xe_x\wedge ye_y$, which is a vector pointing to a direction orthogonal
to the plane as the wedge product is a bit like the cross product. 
Instead of $xe_x$ and $ye_y$ we can use a linear combination of them:
$z=x+iy$ and $\bar{z}=x-iy$. 
In this case the wedge product gets the form
$z\wedge \bar{z}$. Here $z$ is a local coordinate of the complex line. It
is clearer if we denote the local coordinates of the tangent space of the
complex submanifold by $z'_j$ and the local coordinates of the tangent 
space of the mother manifold by $z_m$. Then $z'_j$ is a linear combination
(with complec coefficients) of the local coordinates $z_m$ and a complex
curve expressed as a (1,1)-form in is 
$f(P)dz'\wedge d\bar{z'}$ where $z'=z'_1$
is the only local coordinate of the 1-dimensional submanifold and $P$ is 
the basepoint.
If we take a complex surface as the submanifold, it has two local coordinates
$z'_1$ and $z'_2$ and the form is of the type 
\begin{equation}
f_1(P)f_2(P)dz_1\wedge d\bar{z}_1\wedge dz_2\wedge d\bar{z}_2 
\end{equation}

Consider a (2,2)-form of the type (changing the notation to upper indices)
\begin{equation}
\psi=f(z,\bar{z})dz^1\wedge dz^2\wedge d\bar{z}^3\wedge d\bar{z}^4
\end{equation}
Does such a form exist in the cohomology of $M$ and can such a form 
be obtained as a $\mathbb{Q}$-linear combination of forms of the type (1.3)?

Lemmas 3.1-3.3 are only needed to show that a form of the type (1.4) is
not exact, i.e., zero in the cohomology. 
Lemmas 3.1 and 3.2 show that for
a particular type of a form $\psi$ we get $*\psi$ which is quite similar to
$\psi$ in form. Lemma 3.3 shows that the form is harmonic and thus not exact. 
Actually, we created the form so that both $d\psi=0$ and $d(*\psi)=0$ for 
trivial reasons. There is a special condition in this part of the proof: 
Lemmas 3.1-3.3 needs an assumption 
that the submanifolds are algebraic varieties defined by polynomial 
equations with real coefficients. The important part is that the coefficients 
are real since the proof makes complex conjugation and requires that 
$f(P)^*=f(P^*)$, which is not always true. This assumption means that
the manifold $M$ must be selected in such a way that this assumption holds. 

Lemma 3.4 shows that a form of the type (1.4) exists on $X=K3\times K3$, 
which is a compact complex manifold of complex dimension 4. Lemma 3.3 shows
that the form is not zero in $H^{2,2}(X)$. This is where we need the surface to
be K3 and e.g. an Abelian surface would not work. Complex conjugation sends a 
harmonic form to a harmonic form, and thus we have a wedge of two harmonic 
forms, but the wedge of two harmonic forms is not necessarily harmonic - 
it can be exact. The nowhere vanishing (2,0)-form is needed in order to 
conclude that the wedge product is harmonic. 

Lemma 3.5 proves that 
a form of the type (1.4) cannot arise as a $\mathbb{Q}$-linear combination 
from the cohomology classes of submanifolds. Here we again need a similar
condition as in Lemma 3.1 that $f(P)^*=f(P^*)$. 
This condition follows if the mother manifold and
all submanifolds are algebraic varieties defined by homogeneous polynomial
equations with real coefficients. A restrictive assumption must be added
to the lemma and fulfilling it requires a special form for $M$. 

The manifold $X$ in Lemma 3.4 is not quite what we want. It is a 
submanifold of $\mathbb{P}^3\times \mathbb{P}^3$ but we want an algebraic 
subvariety of $\mathbb{P}^n$ for some $n$. Thus, we replace $X$ by $M$.
No changes are needed to Lemma 3.4.
 
The (2,2)-form from Lemma 3.4 is rational and represents a Hodge class.
This fact seems obvious considering that it is like (1.4) and has only one 
term. However, in 2011 I received a comment from an expert of the field 
stating that the form is not rational. Therefore this issue must be addresses.

Finally I comment the apparent conflict of of the presented result with 
a published Ph.D. thesis [4]. The thesis proves that for some K3xK3 
manifolds the Hodge conjecture holds. The conflict is only apparent since 
the proofs of Lemmas 3.1 and 3.5 require special conditions that 
$f(P)^*=f(P^*)$ and
do not state anything for the general case of K3xK3.

\section {Notations and concepts}

We will use the notations in Kodaira [7] page 147. Local coordinates of a complex manifold $M$ of complex dimension $n$ at the base point $z_0$ are denoted
by $z^1,\dots, z^n$, $\bar{f}(z^1,\dots, z^n)$ is the complex conjugate of 
$f(z^1,\dots,z^n)$, and the Hodge star operation is denoted by $*$.

Let $\varphi$ and $\psi$ be
$C^{\infty}(p,q)$-forms in a complex manifold $M$. 
\begin{equation}
\begin{gathered}
\varphi=\frac {1}{p!q!}\sum \varphi_{\alpha_1,\dots,\alpha_p,\bar{\beta_1},\dots,\bar{\beta_q}}(z)dz^{\alpha_1}\wedge\cdots\wedge dz^{\alpha_p}\wedge d\bar{z}^{\beta_1}\wedge\cdots\wedge d\bar{z}^{\beta_q}\\
\psi=\frac {1}{p!q!}\sum \psi_{\alpha_1,\dots,\alpha_p,\bar{\beta_1},\dots,\bar{\beta_q}}(z)dz^{\alpha_1}\wedge\cdots\wedge dz^{\alpha_p}\wedge d\bar{z}^{\beta_1}\wedge\cdots\wedge d\bar{z}^{\beta_q}\\
\end{gathered}
\end{equation}
The inner product is defined as
\begin{equation}
(\varphi,\psi)(z)=\frac {1}{p!q!}\sum \varphi_{\alpha_1,\dots,\alpha_p,\bar{\beta_1},\dots,\bar{\beta_q}}(z)\bar{\psi}_{\alpha_1,\dots,\alpha_p,\bar{\beta_1},\dots,\bar{\beta_q}}(z)
\end{equation}
and
\begin{equation}
(\varphi,\psi)=\int_M(\varphi,\psi)(z)\frac {\omega^n}{n!}
\end{equation}
where the volume element $\frac {\omega^n}{n!}$ is
\begin{equation}
\frac {\omega^n}{n!}=(-1)^{n(n-1)/2}g(z)dz^1\wedge\cdots\wedge dz^n\wedge d\bar
{z}^1\wedge\cdots\wedge d\bar{z}^n
\end{equation}
Here 
\begin{equation}
g(z)=det(g_{\alpha\bar{\beta}}(z))_{\alpha,\beta=1,\dots,n}
\end{equation}
is given by a Hermitian metric
\begin{equation}
\sum_{\alpha,\beta=1}^n g_{\alpha\bar{\beta}}dz^{\alpha}\otimes d\bar{z}^{\beta}
\end{equation}
and
\begin{equation}
\omega=i\sum_{\alpha,\beta=1}^n g_{\alpha\bar{\beta}}dz^{\alpha}d\bar{z}^{\beta}
\end{equation}
is the associated (1,1)-form, and $\omega^n=\omega\wedge\cdots \wedge \omega$ 
is the $n$-fold product. Then
\begin{equation}
\begin{gathered}
(\varphi,\psi)\frac {\omega^n}{n!}=(i)^n(-1)^{n(n-1)/2}g(z)\frac{1}{p!q!}\sum\varphi_{\alpha_1\dots\alpha_p\bar{\beta}_1\dots\bar{\beta}_q}(z)\\
\times \bar{\psi}_{\alpha_1\dots\alpha_p\bar{\beta}_1\dots\bar{\beta}_q}dz^1\wedge\cdots\wedge dz^n\wedge d\bar{z}^1\wedge\cdots \wedge\bar{z}^n\\
=\frac{1}{p!q!}\sum\varphi_{\alpha_1\dots\alpha_p\bar{\beta}_1\dots\bar{\beta}_q}(z)dz^{\alpha_1}\wedge\cdots\wedge dz^{\alpha_p}\wedge d\bar{z}^{\beta_1}\wedge d\bar{z}^{\beta_q}\wedge *\psi(z)\\
=\varphi(z)\wedge *\psi(z)
\end{gathered}
\end{equation}
where
\begin{equation}
\begin{gathered}
*\psi(z)=(i)^n(-1)^{\frac 12 n(n-1)+(n-p)q}\sum_{A_p,B_q}sgn
\left( \begin{matrix}
A_p & A_{n-p}\cr
B_q & B_{n-q}
\end{matrix} \right)g(z)\bar{\psi}^{A_pB_q}(z)\\
\times dz^{B_{n-q}}\wedge d\bar{z}^{A_{n-p}}
\end{gathered}
\end{equation}
The multi-indices are $A_p=\alpha_1\dots\alpha_p$, $A_{n-p}=\alpha_{p+1}\dots\alpha_n$, $B_q=\beta_1\dots\beta_q$, $B_{n-q}=\beta_{q+1}\dots\beta_n$.
Notice, that Kodaira on page 117 defines the Hodge star differently as
\begin{equation}
(\varphi,\psi)\frac {\omega^n}{n!}=\varphi(z)\wedge *\bar{\psi}(z)
\end{equation}
but the usual definition in the literature is 
\begin{equation}
(\varphi,\psi)vol_x=\varphi(z)\wedge *\psi(z)
\end{equation}
i.e., 
\begin{equation}
(\varphi,\psi)\frac {\omega^n}{n!}=\varphi(z)\wedge *\psi(z)
\end{equation}
as we have defined. 

Let us select local coordinates such that $g_{\alpha\bar{\beta}}(z_0)=\delta_{\alpha\beta}$ where $z_0$ is the base point. Then also $g^{\alpha\bar{\beta}}(z_0)=\delta_{\alpha\beta}$.

We will usually mark indices as upper indices but in a calculation in Lemma 3.5
upper indices get confused with powers and indices are marked as lower indices.

\section{Lemmas and a Theorem}

\begin{lemma}
Let $p=q=2$, $n=4$, and
\begin{equation}
\begin{gathered}
\psi_{12\bar 3\bar 4}(z)=4f(z^1,z^2,\bar{z}^3, \bar{z}^ 4)\\
\psi_{\alpha_1\alpha_2\bar{\beta_1}\bar{\beta_2}}(z)=0\quad if\quad (\alpha_1\alpha_2\bar{\beta_1}\bar{\beta_2})\not=(12\bar{3}\bar{4})
\end{gathered}
\end{equation}
and let $f(z^1,z^2,z^3,z^ 4)$ be a holomorphic function satisfying
$f(z^1,z^2,z^3,z^ 4)^*=f(\bar{z}^1,\bar{z}^2,\bar{z}^3,\bar{z}^4)$.
 Then in local 
coordinates
\begin{equation}
\begin{gathered}
\psi=f(z^1,z^2,\bar{z}^3,\bar{z}^4)dz^1\wedge dz^2 \wedge d\bar{z}^3\wedge d\bar{z}^4\\
*\psi=g(z)f(\bar{z}^1,\bar{z}^2,z^3,z^4)dz^3\wedge dz^4 \wedge d\bar{z}^1\wedge d\bar{z}^2\\
\end{gathered}
\end{equation}
\end{lemma}

\begin{proof}
The first claim is obvious since
\begin{equation}
\psi=\frac 14 \psi_{12\bar{3}\bar{4}}(z)dz^1\wedge dz^2 \wedge d\bar{z}^3\wedge d\bar{z}^4\\
\end{equation}
For the second claim we calculate
\begin{equation}
\begin{gathered}
*\psi=(i)^4(-1)^{\frac 12 4\cdot 3+2\cdot 2}sgn
\left( \begin{matrix}
1 & 2 & 3 & 4\cr
3 & 4 & 1 & 2
\end{matrix} \right)
g(z)\psi^{12\bar{3}\bar{4}}(z)dz^3\wedge dz^4\wedge d\bar{z}^1\wedge d\bar{z}^2
\end{gathered}
\end{equation}
Since 
\begin{equation}
\begin{gathered}
\psi^{\alpha_1\alpha_2\bar{\beta_1}\bar{\beta_2}}(z)=\sum_{\lambda_1,\lambda_2,\mu_1,\mu_2}g^{\bar{\lambda_1}\alpha_1}g^{\bar{\lambda_2}\alpha_2}g^{\bar{\beta_1}\mu_1}g^{\bar{\beta_2}\mu_2}\bar{\psi}_{\lambda_1\lambda_2\bar{\mu_1}\bar{\mu_2}}(z)\\
=\sum_{\lambda_1,\lambda_2,\mu_1,\mu_2}\delta_{\lambda_1\alpha_1}\delta_{\lambda_2\alpha_2}\delta_{\beta_1\mu_1}\delta_{\beta_2\mu_2}\bar{\psi}_{\lambda_1\lambda_2\bar{\mu_1}\bar{\mu_2}}(z)\\
=\psi_{\alpha_1\alpha_2\bar{\beta_1}\bar{\beta_2}}(z)=\left\{\aligned
&\bar{\psi}_{12\bar{3}\bar{4}}(z)\quad {if}\quad (\alpha_1\alpha_2\bar{\beta_1}\bar{\beta_2})=(12\bar{3}\bar{4})\\
&0\quad {otherwise}\endaligned\right.
\end{gathered}
\end{equation}
Thus
\begin{equation}
\begin{gathered}
*\psi=g(z)\bar{\psi}_{12\bar{3}\bar{4}}(z)dz^3\wedge dz^4\wedge d\bar{z}^1\wedge \bar{z}^2\\
=g(z)\bar{f}(z^1,z^2,\bar{z}^3,\bar{z}^4)dz^3\wedge dz^4\wedge d\bar{z}^1\wedge \bar{z}^2
\end{gathered}
\end{equation}
and the assumption on $f(z^1,z^2,z^3,z^4)$ gives
\begin{equation} 
\bar{f}(z^1,z^2,\bar{z}^3,\bar{z}^4)=f(z^1,z^2,\bar{z}^3,\bar{z}^4)^*
=f(\bar{z}^1,\bar{z}^2,z^3,z^4)
\end{equation} 
\end{proof}

\begin{lemma}
Let $p=q=2$, $n=4$, and
\begin{equation}
\begin{gathered}
\psi_{12\bar 3\bar 4}(z)=4f(z^1,z^2,\bar{z}^3, \bar{z}^ 4)\\
\psi_{\alpha_1\alpha_2\bar{\beta_1}\bar{\beta_2}}(z)=0\quad if\quad (\alpha_1\alpha_2\bar{\beta_1}\bar{\beta_2})\not=(12\bar{3}\bar{4})
\end{gathered}
\end{equation}
and let $f(z^1,z^2,z^3,z^4)$ be a holomorphic function
satisfying
$f(z^1,z^2,z^3,z^ 4)^*=f(\bar{z}^1,\bar{z}^2,\bar{z}^3,\bar{z}^4)$.
We can select
the metric such that in local 
coordinates
\begin{equation}
\begin{gathered}
\psi=f(z^1,z^2,\bar{z}^3,\bar{z}^4)dz^1\wedge dz^2 \wedge d\bar{z}^3\wedge d\bar{z}^4\\
*\psi=f(\bar{z}^1,\bar{z}^2,z^3,z^4)dz^3\wedge dz^4 \wedge d\bar{z}^1\wedge d\bar{z}^2\\
\end{gathered}
\end{equation}
\end{lemma}

\begin{proof}
At the base point $z_0$ we can select the metric 
$g(z_0)=\det (g_{\alpha\bar{\beta}})_{\alpha,\beta=1,\dots,n}$ 
such that  $g_{\alpha\bar{\beta}}=\delta_{\alpha\beta}$. Then 
$g(z_0)=\det (g_{\alpha\bar{\beta}})_{\alpha,\beta=1,\dots,n}=1$.
We can make the same selection at all points $z$ and thus $g(z)=1$.
\end{proof}

\begin{lemma}
Let $\psi$ be a (2,2)-form in a complex manifold $M$. Let $\psi$ be expressed
in local coordinates as 
\begin{equation}
\psi=f(z^1,z^2,\bar{z}^3,\bar{z}^4)dz^1\wedge dz^2 \wedge d\bar{z}^3\wedge d\bar{z}^4
\end{equation}
where $f(z^1,z^2,\bar{z}^3,\bar{z}^4)$ is holomorphic
satisfying
$f(z^1,z^2,z^3,z^ 4)^*=f(\bar{z}^1,\bar{z}^2,\bar{z}^3,\bar{z}^4)$.
Let us assume that the
complex dimension of $M$ is four, $M$ is compact and without boundary (i.e., 
closed). Then
$\bigtriangleup \psi=0$. 
\end{lemma}

\begin{proof}
For compact manifolds without boundary 
\begin{equation}
\bigtriangleup \psi=0 \quad \Leftrightarrow \quad d\psi=0\quad {and}\quad \delta \psi=0 
\end{equation}
Let us calculate $d\psi$
\begin{equation}
d\psi=\sum_{j=1}^4\frac{\partial f}{\partial z^j}(z^1,z^2,\bar{z}^3,\bar{z}^4)dz^j\wedge dz^1\wedge dz^2 \wedge d\bar{z}^3 \wedge d\bar{z}^4=0
\end{equation}
since $dz^j\wedge dz^1=0$ if $j=1$, $dz^j\wedge dz^2=0$ if $j=2$, and
\begin{equation}
\begin{gathered}
\frac{\partial f}{\partial z^3}(z^1,z^2,\bar{z}^3,\bar{z}^4)=0\\
\frac{\partial f}{\partial z^4}(z^1,z^2,\bar{z}^3,\bar{z}^4)=0
\end{gathered}
\end{equation}
since in the coordinate system $(z^1,\dots,z^n,\bar{z}^1,\dots ,\bar{z}^n)$ the
coordinates 
$z^j$ and $\bar{z}^j$ are considered independent. The second assertation follows
from the definition of the codifferential: if $*:\Omega^k\to \Omega^{n-k}$
is the Hodge star operator then $\delta:\Omega^k\to \Omega^{k-1}$ is defined
by
\begin{equation}
\delta\psi=(-1)^{n(k+1)+1}(*d*)\psi
\end{equation}
Inserting $n=4$, $k=2$, yields
\begin{equation}
\delta\psi=-(*d*)\psi=-*d(*\psi)
\end{equation}
We may assume that the metric is chosen such that $g(z)=1$. By Lemma 3.2 
\begin{equation}
*\psi=f(\bar{z}^1,\bar{z}^2,z^3,z^4)dz^3\wedge dz^4 \wedge d\bar{z}^1\wedge d\bar{z}^2
\end{equation}
As in the previous case, we conclude that
\begin{equation}
d(*\psi)=\sum_{j=1}^4\frac{\partial f}{\partial z^j}(\bar{z}^1,\bar{z}^2,z^3,z^4)dz^j\wedge dz^3\wedge dz^4 \wedge d\bar{z}^1 \wedge d\bar{z}^2=0
\end{equation}

\end{proof}

\begin{lemma}
Let $X$ be a product of two K3 surfaces. Then $X$ allows a (2,2)-form of the 
type
\begin{equation}
f(z^1,z^2)f(\bar{z}^3,\bar{z}^4)dz^1\wedge dz^2\wedge d\bar{z}^3\wedge d\bar{z}^4
\end{equation}
where $f(z^1,z^2)$ is holomorphic and nowhere vanishing.
\end{lemma}

\begin{proof}
The existence of a nowhere vanishing 2-form is often taken as the definition 
of a K3 surface, the additional condition guaranteeing that a 
2-dimensional complex manifold $X$ is a K3 surface is that the manifold $X$
is connected. This nowhere vanishing 2-form $\lambda$ is the generator of 
$H^{2,0}(X)$, i.e., every other element $\alpha\in H^{2,0}(X)$ can be 
expressed as $\alpha=c\lambda$ for some $c\in \mathbb{C}$. The complex 
conjugate $\bar{\lambda}$ of $\lambda$ is the generator of $H^{0,2}(X)$
as is shown by the Hodge duality pairing. Let the (2,0)-form $\lambda$ be 
expressed in local coordinates as
\begin{equation}
\lambda=f(z^1,z^2)dz^1\wedge dz^2
\end{equation}
Then $f(z^1,z^2)$ is holomorphic and nowhere vanishing. 
If $M=K3\times K3$ then there are two nowhere vanishing 2-forms 
$\lambda_1$ and $\lambda_2$. Let us remember that a complex K3-surface is
compact and as a real manifold it is a 4-dimensional closed manifold.
Thus, $M$ is closed and we can use Lemma 3.3.  
We can make a (2,2)-form as a wedge product
of 
\begin{equation}
\lambda_1=f(z^1,z^2)dz^1\wedge dz^2\\
\end{equation}
and
\begin{equation}
\bar{\lambda}_2=\bar{f}(z^3,z^4)d\bar{z}^3\wedge d\bar{z}^4
=f(\bar{z}^3,\bar{z}^4)d\bar{z}^3\wedge d\bar{z}^4
\end{equation}
as $f(z^3,z^4)$ is holomorphic. Thus, we define
\begin{equation}
\begin{gathered}
\begin{aligned}
\psi&=f(z^1,z^2)dz^1\wedge dz^2\wedge f(\bar{z}^3,\bar{z}^4)d\bar{z}^3\wedge d\bar{z}^4\\
&=f(z^1,z^2)f(\bar{z}^3,\bar{z}^4)dz^1\wedge dz^2\wedge d\bar{z}^3\wedge d\bar{z}^4
\end{aligned}
\end{gathered}
\end{equation}  
The form $\psi$ satisfies $\bigtriangleup \psi=0$ by Lemma 3.3.
We still have to show that it is not zero in $H^{2,2}(X)$.
Typical ways to show this are calculating periods or intersection products,
but we will show it differently. 
For any complex manifold $M$
the Hodge star gives an 
isomorphism from $H^k(M)$ to $H^{n-k}(M)$ and the Poincar\'e duality gives 
an isomprphism from $H_k(M)$ to $H^{n-k}(M)$. Especially, when $n=4$, 
$H^2(X)$ and $H_2(X)$ are isomorphic and if $X=S_1\times S_2$ 
$H^2(S_i)$ and $H_2(S_i)$, $i=1,2$ are isomorphic (for any coefficients, so
coefficients are supressed in the notations). 
Let this isomorphism be
$\Psi:H_2(X)\to H^2(X)$.
Let $C$ be a 2-chain in $S_1$ 
and let $(C,p_t)\in X$ be a family of pairs where $p_t$ is a path in 
$S_2$. This family $C\times [0,1]$ defines a homotopy from $(C,t_0)$ to
$(C,t_1)$. By the homomorphism we have also a family 
$\lambda_1\wedge \bar\lambda_2(P_t)$ parametrized by $t\in [0,1]$. Then $P_t$
defines a path in $S_2$. This is where we need the nowhere vanishing
(2,0)-form. If there existed a point $P_1$ such that 
$\bar\lambda_2(P_1)=0$ then $\psi=0$ at any point $(Q,P_1)\in S_1\times S_2$.
As the preimage of zero is a point, this would yield a homotopy from $C$ to 
a point, i.e., $(C,p_0)=0$ in $\pi_2(X)$. Consequently $(C,p_0)=0$ in $H_2(X)$
and $\Psi(C,p_0)=(\lambda_1,P_0)=0$ in $H^2(X)$. 
As there is no such point $P_1$,
$(C,p_0)$ is not contractible in $X$ and is therefore nonzero in $H_2(X)$.   
Its image under $\Psi$ is therefore nonzero in $H^2(X)$.
Thus, $\psi$ is a harmonic form in $H^{2,2}(X)$.
\end{proof}

\begin{lemma}
A (2,2)-form of the type $f(z,\bar{z})dz^1\wedge dz^2\wedge d\bar{z}^3\wedge d\bar{z}^4$ is not a linear combination with rational coefficients
of cohomology classes deriving from complex submanifolds of complex 
codimension 2 in a complex submanifold of complex dimension 4 assuming two
conditions:
1) the mother manifold is an algebraic variety defined by homogeneous 
polynomial equations with real coefficients, and 2) 
all submanifolds are algebraic varieties 
defined by homogeneous polynomial equations with real coefficient.
\end{lemma}

\begin{proof}
Let $N$ be a complex manifold of dimension $n$ and let $M$ be a complex 
submanifold of $N$ of complex dimension $k$. Let $i:M\to N$ be the inclusion
map and the complex structure of $M$ be induced by the complex structure of
$N$. Let $P\in M$ be a point in $M$ and $(x'^1,y'^1,\dots,x'^k,y'^k)$ be 
a local
coordinate system in $TM_P$ where $M$ is considered as a real manifold of 
dimension $2k$. The coordinate system can be completed to a coordinate system
of $TN_{i(P)}=TN_P$ by adding $2(n-k)$ coordinate vectors. This yields a 
coordinate system 
\begin{gather*}
(x',y')=(x'^1,y'^1,\dots,x'^k,y'^k,x'^{k+1},y'^{k+1},\dots, x'^n,y'^n)
\end{gather*} 
to $TN_P$ where $P=(0,0)=(0,\dots,0)$. 
We can assume that the coordinates are orthonormal. 
The manifold $N$ can be considered as a real manifold of dimension $2n$ and
$TN_P$ be given a local coordinate system  
\begin{gather*}
(x,y)=(x^1,y^1,\dots, x^n,y^n)\quad,\quad P=i(P)=(0,0)
\end{gather*} 
There is 
an orthonormal coordinate transform $A:\mathbb{R}^{2n}\to \mathbb{R}^{2n}$,
$A\in SO(2n)$, such that
\begin{equation}
[x' \ y']^T=A[x \ y]^T
\end{equation}
The coordinate systems at $TN_P$ can be chosen such that
\begin{equation}
\begin{gathered}
\begin{aligned}
&z^j=x^j+iy^j\quad, \quad \bar{z}^j=x^j-iy^j\quad j=1,\dots, n\\
&z'^j=x'^j+iy'^j\quad, \quad \bar{z}'^j=x'^j-iy'^j\quad j=1,\dots, n
\end{aligned}
\end{gathered}
\end{equation}
and there is an orthonormal coordinate transform $B:\mathbb{C}^n\to \mathbb{C}^n$ such that for
\begin{equation}
\begin{gathered}
(z,\bar{z})=(z^1,\bar{z}^1,\dots,z^n,\bar{z}^n)\\
(z,\bar{z}')=(z'^1,\bar{z}'^1,\dots,z'^n,\bar{z}'^n)
\end{gathered}
\end{equation}
the transform takes $(z,\bar{z})$ to $(z',\bar{z}')$, 
$[z' \ \bar{z}']^T=B[z \ \bar{z}]^T$. As the complex structure of $M$ is 
induced from the complex structure of $N$, the following statement holds
\begin{equation}
\begin{gathered}
{if}\quad (z',\bar{z}')=(x'^1,x'^1,\dots, x'^n,x'^n)\quad i.e. \quad y'^j=0 \ {for} \ j=1,\dots, n\\
{then}\quad (z,\bar{z})=(x^1,x^1,\dots, x^n,x^n)\quad i.e. \quad y^j=0 \ {for} \ j=1,\dots, n
\end{gathered}
\end{equation}
It follows that if
\begin{equation}
z'^j=\sum_m a_{m,j}z^m
\end{equation}
then 
\begin{equation}
\bar{z}'^j=\sum_m a_{m,j}\bar{z}^m
\end{equation}
that is, $a_{m,j}\in \mathbb{R}$ for all $j$ and $m$. The Poincar\'e dual $[M]$
of $M$ considered as a $2k$-dimensional real manifold in the $2n$-dimensional
real manifold $N$ satisfies the condition that $[M]$ capped with the 
fundamental class of $N$ is the homology class of $M$. Thus, $[M]$ is a form
of the type
\begin{equation}
\varphi(x',y')dx'^{k+1}\wedge dy'^{k+1}\wedge \cdots \wedge dx'^n\wedge dy'^n
\end{equation}
Since $x'^j$ and $y'^j$ are independent, the class $[M]$ is represented at
$p$ by a form of the type
\begin{equation}
\varphi(z',\bar{z}')dz'^{k+1}\wedge d\bar{z}'^{k+1}\wedge \cdots \wedge dz'^n\wedge d\bar{z}'^n
\end{equation}
In coordinates $\{ z^1,\bar{z}^1,\dots,z^n,\bar{z}^n\}$ the form can be
expressed by inserting $dz'^j$ and $d\bar{z}'^j$ as linear combinations 
(with real coefficients) of $dz^k$ and $d\bar{z}^k$, $k=1,\dots, n$. 

Let us now take a two dimensional complex plane in a four dimensional
complex manifold. We will use lower indices in this calculation as
we obtain squares at some places 
and upper indices get easily confused with powers. 
The two local coordinates $z''_1,z''_2$ of a complex plane in a 
four dimensional complex manifold can be expressed in
the local coordinates $z_1,z_2,z_3,z_4$ of the mother manifold as
\begin{equation}
\begin{gathered}
z''_1=b_{11}z_1+b_{12}z_2+b_{13}z_3+b_{14}z_4\\
z''_1=b_{21}z_1+b_{22}z_2+b_{23}z_3+b_{24}z_4
\end{gathered}
\end{equation}
where $b_{ij}$ are complex numbers. 
We want to see what terms come from 
\begin{equation}
\begin{gathered}
dz''_1\wedge d\bar{z''}_1\wedge dz'_2\wedge d\bar{z'}_2
\end{gathered}
\end{equation}
when the local coordinates $z''_i$ of the submanifold are replaced by
the local coordinates $z_i$ of the mother manifold.
We get 36 terms the form
$dz_{m_1}\wedge dz_{m_2}\wedge d\bar {z}_{m_3}\wedge d\bar {z}_{m_4}$
where $m_1<m_2$ and $m_3<m_4$ because if f $m_1=1$, $m_1$ can be $2,3,4$, 
if $m_1=2$, $m_2$ can be $3,4$, and if $m_1=3$ then $m_2=4$. 
This means that there are six possibilities
for $m_1,m_2$. There are also six possibilities for $m_3,m_4$. Together
there are 36 possibilities and thus 36 different terms. 
While the terms can be calculated we get simpler expressions 
by making a linear transform of the
local coordinates $z''_1, z''_2$
\begin{equation}
\begin{gathered}
z'_1=\alpha(-b_{22}z''_1+b{12}z''_2)\\
z'_2=\alpha(b_{21}z''_1-b{11}z''_2)\\
\alpha=(b_{21}b_{12}-b_{11}b_{22})^{-1}
\end{gathered}
\end{equation}
Then
\begin{equation}
\begin{gathered}
z'_1=z_1+a_{13}z_3+a_{14}z_4\\
z'_2=z_2+a_{23}z_3+a_{24}z_4\\
a_{13}=\alpha(b_{12}b_{23}-b_{13}b_{22})\quad,\quad
a_{14}=\alpha(b_{12}b_{24}-b_{14}b_{22})\\
a_{23}=\alpha(b_{21}b_{13}-b_{11}b_{23})\quad,\quad
a_{24}=\alpha(b_{21}b_{14}-b_{11}b_{24})\\
\end{gathered}
\end{equation}
These new coordinates are not orthogonal, but that does not matter
here: we just want to see what terms come from 
\begin{equation}
\begin{gathered}
dz'_1\wedge d\bar{z'}_1\wedge dz'_2\wedge d\bar{z'}_2
\end{gathered}
\end{equation}
when the $z_i$ coordinates are inserted. Again there are 36 terms, but 
the expressions are shorter than for the original coordinates $z''_i$. 
The first six terms are:
\begin{equation}
\begin{gathered}
=-dz_1\wedge dz_2\wedge d\bar{z}_1\wedge d\bar{z}_2
-|a_{13}|^2 dz_2\wedge dz_3\wedge d\bar{z}_2\wedge d\bar{z}_3
-|a_{14}|^2 dz_2\wedge dz_4\wedge d\bar{z}_2\wedge d\bar{z}_4\\
-|a_{23}|^2 dz_1\wedge dz_3\wedge d\bar{z}_1\wedge d\bar{z}_3
-|a_{24}|^2 dz_1\wedge dz_4\wedge d\bar{z}_1\wedge d\bar{z}_4\\
-|a_{23}a_{14}-a_{24}a_{13}|^2 dz_3\wedge dz_4\wedge d\bar{z}_3\wedge d\bar{z}_4
\end{gathered}
\end{equation}
Terms 7-16 numbered in this order:
\begin{equation}
\begin{gathered}
+a^*_{13}dz_1\wedge dz_2\wedge d\bar{z}_2\wedge d\bar{z}_3
+a^*_{14}dz_1\wedge dz_2\wedge d\bar{z}_2\wedge d\bar{z}_4\\
+a_{13}dz_2\wedge dz_3\wedge d\bar{z}_1\wedge d\bar{z}_2
+a_{14}dz_2\wedge dz_4\wedge d\bar{z}_1\wedge d\bar{z}_2\\
-a_{13}a^*_{14} dz_2\wedge dz_3\wedge d\bar{z}_2\wedge d\bar{z}_4
-a^*_{13}a_{14} dz_2\wedge dz_4\wedge d\bar{z}_2\wedge d\bar{z}_3\\
+a_{23}(a^*_{23}a^*_{14}-a^*_{24}a^*_{13}) dz_1\wedge dz_3\wedge d\bar{z}_3\wedge d\bar{z}_4\\
+a_{24}(a^*_{23}a^*_{14}-a^*_{24}a*_{13}) dz_1\wedge dz_4\wedge d\bar{z}_3\wedge d\bar{z}_4\\
+a^*_{23}(a_{23}a_{14}-a_{24}a_{13}) dz_3\wedge dz_4\wedge d\bar{z}_1\wedge d\bar{z}_3\\
+a^*_{24}(a_{23}a_{14}-a_{24}a_{13}) dz_3\wedge dz_4\wedge d\bar{z}_1\wedge d\bar{z}_4
\end{gathered}
\end{equation}
Terms 17-22:
\begin{equation}
\begin{gathered}
-a_{24}a^*_{23}dz_1\wedge dz_4\wedge d\bar{z}_1\wedge d\bar{z}_3
-a^*_{24}a_{23}dz_1\wedge dz_3\wedge d\bar{z}_1\wedge d\bar{z}_4\\
-a_{13}(a^*_{23}a^*_{14}-a^*_{24}a^*_{13}) dz_2\wedge dz_3\wedge d\bar{z}_3\wedge d\bar{z}_4\\
-a_{14}(a^*_{23}a^*_{14}-a^*_{24}a^*_{13}) dz_2\wedge dz_4\wedge d\bar{z}_3\wedge d\bar{z}_4\\
-a^*_{13}(a_{23}a_{14}-a_{24}a_{13}) dz_3\wedge dz_4\wedge d\bar{z}_2\wedge d\bar{z}_3\\
-a^*_{14}(a_{23}a_{14}-a_{24}a_{13}) dz_3\wedge dz_4\wedge d\bar{z}_2\wedge d\bar{z}_4
\end{gathered}
\end{equation}
Terms 23-34:
\begin{equation}
\begin{gathered}
+a^*_{23}a_{13}dz_2\wedge dz_3\wedge d\bar{z}_1\wedge d\bar{z}_3
+a^*_{23}a_{14}dz_2\wedge dz_4\wedge d\bar{z}_1\wedge d\bar{z}_3\\
+a^*_{24}a_{13}dz_2\wedge dz_3\wedge d\bar{z}_1\wedge d\bar{z}_4
+a^*_{24}a_{14}dz_2\wedge dz_4\wedge d\bar{z}_1\wedge d\bar{z}_4\\
+a_{23}a^*_{13}dz_1\wedge dz_3\wedge d\bar{z}_2\wedge d\bar{z}_3
+a_{23}a^*_{14}dz_1\wedge dz_3\wedge d\bar{z}_2\wedge d\bar{z}_4\\
+a_{24}a^*_{13}dz_1\wedge dz_4\wedge d\bar{z}_2\wedge d\bar{z}_3
+a_{24}a^*_{14}dz_1\wedge dz_4\wedge d\bar{z}_2\wedge d\bar{z}_4\\
-a^*_{23} dz_1\wedge dz_2\wedge d\bar{z}_1\wedge d\bar{z}_3
-a^*_{24} dz_1\wedge dz_2\wedge d\bar{z}_1\wedge d\bar{z}_3\\
-a_{23} dz_1\wedge dz_3\wedge d\bar{z}_1\wedge d\bar{z}_2
-a_{24} dz_1\wedge dz_4\wedge d\bar{z}_1\wedge d\bar{z}_2
\end{gathered}
\end{equation}
We are mainly interested in these two last terms, terms 35 and 36:
\begin{equation}
\begin{gathered}
+(a_{23}a_{14}-a_{24}a_{13})dz_3\wedge dz_4\wedge d\bar{z}_1\wedge d\bar{z}_2\\
+(a^*_{23}a^*_{14}-a^*_{24}a^*_{13})dz_1\wedge dz_2\wedge d\bar{z}_3\wedge d\bar{z}_4
\end{gathered}
\end{equation}

We assume that 
the mother manifold is an algebraic variety defined by a finite number of
homogeneus polynomial
equations of a finite number of complex parameters $s_i$ and real 
coefficients. An example of such a variety is the Fermat quadratic 
surface defined by the homogeneous polynomial equation
\begin{equation}
\begin{gathered}
s_1^4+s_2^4+s_3^4+s_4^4=0
\end{gathered}
\end{equation}
In thie polynomial equation the coefficients are all integers (all are 1s).

If $P_o=(s_{10},s_{20},s_{30},s_{40})$ is a chosen basepoint, the local
coordinates of the manifold can be chosen as $z_i=s_i-s_{i0}$. 
Let us consider a (2,2)-form
\begin{equation}
\begin{gathered}
f(P)dz_1\wedge dz_2\wedge d\bar{z}_3\wedge d\bar{z}_4 
\end{gathered}
\end{equation}
that is linear combination with rational coefficients of classes of 
$n$ complex submanifolds. Each submanifold gives 36 terms $i$ of the type
\begin{equation}
\begin{gathered}
\varphi_j(P)\varphi_k(P) A_{j,k,i} dz_{m_1}\wedge dz_{m_2}\wedge d\bar{z}_{m_3}\wedge d\bar{z}_{m_4} 
\end{gathered}
\end{equation}
Here $A_{j,k,i}$ is the coefficient of 
$dz_{m_1}\wedge dz_{m_2}\wedge d\bar{z}_{m_3}\wedge d\bar{z}_{m_4}$
for two (1,1)-forms $j$ and $k$ and $i$ is numbered as in (3.36-3.40)
where we have explicitly written the coefficients $A_{1,2,i}$, $i=1,\dots, 36$.Thus,
for $j=1$ and $k=2$ the coefficients of the last two terms are
\begin{equation}
\begin{gathered}
A_{1,2,35}=a_{23}a_{14}-a_{24}a_{13}\\
A_{1,2,36}=a^*_{23}a^*_{14}-a^*_{24}a^*_{13}
\end{gathered}
\end{equation}
where $a_{mn}$ are calculated for the submanifold corresponding to two
(1,1)-forms that we numbered as $j=$ and $k=2$. For different values
of $j$ and $k$ we get different $a_{mn}$ but it was inconvenient to
add the indices $j,k$ to (3.36)-(3.40).
 
If we use the coordinates $z''_i$ instead of the simpler coordinates
$z'_i$ and calculate $A_{1,2,35}$ and $A_{1,2,36}$, they are
\begin{equation}
\begin{gathered}
A_{1,2,35}=-b_{13}b_{24}b^*_{11}b^*_{22}+b_{14}b_{23}b^*_{11}b^*_{22}
+b_{13}b_{24}b^*_{12}b^*_{21}-b_{14}b_{23}b^*_{12}b^*_{21}\\
A_{1,2,36}=-b^*_{13}b^*_{24}b_{11}b_{22}+b^*_{14}b^*_{23}b_{11}b_{22}
+b*_{13}b^*_{24}b_{12}b_{21}-b^*_{14}b^*_{23}b_{12}b_{21}
\end{gathered}
\end{equation}
Clearly, always holds
\begin{equation}
\begin{gathered}
A_{j,k,36}={A_{j,k,35}}^*
\end{gathered}
\end{equation} 
In order for a (2,2)-form of the type (3.42) 
to be a $\mathbb{Q}$-linear combination forms
corresponding to submanifolds (of this type) in each point
$P$ must hold 35 equations of the type
\begin{equation}
\begin{gathered}
\sum_{j,k}\varphi_j(P)\varphi_j(P)A_{j,k,i}=0 \quad i=1,\dots,35
\end{gathered}
\end{equation}  
and the last sum must be non-zero
\begin{equation}
\begin{gathered}
\sum_{j,k}c_{jk}\varphi_j(P)\varphi_j(P)A_{j,k,36}\not=0
\end{gathered}
\end{equation}  
where $c_{jk}$ are rational coefficients.
This is a set of linear equations for the unknowns $c_{jk}$.
It is clear that in a chosen basepoint $P_0$ we can solve 
these linear equations and find complex numbers $c_{jk}$
that satisfy these 36 equations, provided that there are at least 36 
independent submanifolds. It is not clear if we can find rational
numbers $c_{jk}$ that satisfy all equations, but let
us assume that we have found rational
$c_{jk}$ that satisfy these equations at the point $P_0$. The question
is if these same $c_{jk}$ can satisfy these 36 equations in every other
point $P$.   

It is easily shown that they cannot, assuming what we have assumed that
the homogeneous polynomial equations defining the mother manifold and 
the submanifolds have real coefficients. With real coefficients we can
find another point $P_1$ as a complex conjugate of $P_0$. 

Thus, let us take the second point as 
$P_1=(s^*_{10},s^*_{20},s^*_{30},s^*_{30})$. Then $P_1=P^*_0$. 
The mother manifold and submanifolds are assumed to be algebraic varieties
defined by a finite number of homogeneous polynomials with real coefficients.
Therefore if $P_0$ is a solution
to the homogeneous polynomial equations, then so it $P^*_0$. 
The functions $\varphi_j(P)$ have the same form at $P=P_0$ and $P=P_1$
but the values differ: $\varphi_j(P_1)=\varphi_j(P_0^*)=\varphi_j(P_0)^*$ 
because the coefficients of the homogeneous polynomial equations
are real. This is the essence: in order to show that the 36 equations
cannot be satisfied at every point $P$ we need to find another point $P_1$
which is far from $P_0$ (in the close vicinity of $P_0$ the 36 equations are
satisfied) and be able to calculate what the 36 equations are at the 
point $P-1$. It is not likely that the 36 equations could be satisfied
at each point $P$ even if we would not make the restrictive assumption,
but it is difficult to find a point $P_1$ without making this assumption.

There are two ways to transfer the local coordinates from $P_0$ to $P_1$.
One way is to choose similar local coordinates $z_i$, $i=1,\dots,4$, 
and $z''_i$, $j=1,2$, in the basepoint $P_1$ as in $P_0$. In this case
the form $dz_1\wedge dz_2\wedge d\bar{z}_3\wedge d\bar{z}_4$ does not
change in the transformation and the
numbers $A_{jk,i}$, $i=1,\dots, 36$, do not change:
\begin{equation}
\begin{gathered}
A_{j,k,36}(P_1)=A_{j,k,36}(P_0)
\end{gathered}
\end{equation}  
where $A_{j,k,i}(P)$ means $A_{j,k,i}$ calculated at the point $P$.
The other method is to take conjugates. In this case both 
the form $dz_1\wedge dz_2\wedge d\bar{z}_3\wedge d\bar{z}_4$ and the
numbers $A_{jk,i}$, $i=1,\dots, 36$ change to complex conjugates.

Let us follow the first way. 
The equation 36 at the point $P=P_1$ 
gives the coefficient of the term 
$dz_1\wedge dz_2\wedge d\bar{z}_3\wedge d\bar{z}_4$. The coefficient is
\begin{equation}
\begin{gathered}
\sum_{j,k}c_{jk}\varphi_j(P_1)\varphi_j(P_1)A_{j,k,36}(P_1)
=\sum_{j,k}c_{jk}\varphi_j(P_0)^*\varphi_j(P_0)^*A_{j,k,36}(P_0)\\
=\left(\sum_{j,k}c_{jk}\varphi_j(P_0)\varphi_j(P_0)\left(A_{j,k,36}(P_0)\right)^*\right)^*\\
=\left(\sum_{j,k}c_{jk}\varphi_j(P_0)\varphi_j(P_0)A_{j,k,35}(P_0)\right)^*=0
\end{gathered}
\end{equation}
That is, at the point $P_1$ the form 
$dz_1\wedge dz_2\wedge d\bar{z}_3\wedge d\bar{z}_4$ vanishes.
The coefficient that does not vanish at $P_1$ is the coefficient of the term 
$dz_3\wedge dz_4\wedge d\bar{z}_1\wedge d\bar{z}_2$.

The second way corresponds to finding the form at $P_1$ by conjugating
the form at $P_0$. In this case
\begin{equation}
\begin{gathered}
\left(\sum_{j,k}c_{jk}\varphi_j(P_0)\varphi_j(P_0)A_{j,k,36}(P_0)dz_1\wedge dz_2\wedge d\bar{z}_3\wedge d\bar{z}_4\right)^*\\
=\sum_{j,k}c_{jk}\varphi_j(P_1)\varphi_j(P_1)A_{j,k,36}(P_1)\left(dz_1\wedge dz_2\wedge d\bar{z}_3\wedge d\bar{z}_4\right)^*\\
=\sum_{j,k}c_{jk}\varphi_j(P_1)\varphi_j(P_1)A_{j,k,36}(P_1)d\bar {z}_1\wedge d\bar{z}_2\wedge dz_3\wedge dz_4\\
=\sum_{j,k}c_{jk}\varphi_j(P_1)\varphi_j(P_1)A_{j,k,36}(P_1)dz_3\wedge dz_4\wedge d\bar {z}_1\wedge d\bar{z}_2
\end{gathered}
\end{equation}
In this method the coefficient of the form stays as nonzero, but it is not the
coeffient of $dz_1\wedge dz_2\wedge d\bar{z}_3\wedge d\bar{z}_4$ at $P_1$.
It is the coefficient of $dz_3\wedge dz_4\wedge d\bar {z}_1\wedge d\bar{z}_2$.
The coeffient of $dz_1\wedge dz_2\wedge d\bar{z}_3\wedge d\bar{z}_4$ at $P_1$
is zero also in this method. Naturally, both ways give the same result. 

We conclude that the form (3.42) cannot be created as a linear combination
of terms (3.43). There should be at least two nonzero terms, those 
corresponding to $A_{j,k,36}$ and to $A_{j,k,35}$.
\end{proof}

\begin{theorem}
The algebraic variety
\begin{equation}
\begin{gathered}
M=\{(s_1:s_2:s_3:s_4:s_5:s_6:s_7:s_8)|s_1^4+s_2^4+s_3^4+s_4^4=0,\\
s_5^4+s_6^4+s_7^4+s_8^4=0,s_4^4+s_5^4=0\}
\end{gathered}
\end{equation}
is an algebraic subvariety of the complex projective space $\mathbb{P}^7$ of 
codimension 3 and a complex manifold of dimension 4. 
There is a Hodge class in $H^{2,2}$ that cannot be represented as a 
$Q$-linear combination of the classes of algebraic subvarieties of $M$.
\end{theorem}

\begin{proof}
Clearly $M$ is a submanifold of $\mathbb{P}^7$ and it is an algebraic variety.
It contains two copies of the Fermat quadratic surface, which is a $K_3$ 
surface. 
The submanifolds of dimension 1 for the Fermat quadratic surface are
defined by adding one homogeneous polynomial equation. The dimension
of the space of (1,1)-forms in $K_3$ is 20, thus if we find 20 polynomial
equations that give an essentially different submanifold, we have represented
all (1,1)-forms by an algebraic subvarity.  
By inspectation, the homogeneous polynomial equations are
\begin{equation}
\begin{gathered}
s_i^4=0 \quad i=1,2,3,4\\
s_i^4=1 \quad i=1,2,3,4\\
s_i^4+s_j^4=0 \quad (i,j)=(1,2),(1,3),(1,4),(2,3),(2,4),(3,4)\\
s_i^4+s_j^4=1 \quad (i,j)=(1,2),(1,3),(1,4),(2,3),(2,4),(3,4)
\end{gathered}
\end{equation}
Clearly, fixing one of $s_i$ to zero gives a submanifold and the submanifold
is different for each $i=1,2,3,4$. In the complex projective space 
$\mathbb{P}^1$ there are only two numbers, zero and one. Thus, setting
$s_i=1$ gives another set of four different submanifolds and there are
no more equations tying only one $s_i$. Setting the sum of two $s_i$ 
variable terms together
gives six different submanifolds if set the sum to zero and another six
if we set the sum to one. There could be more equations of this type, but
as these 20 equations give 20 different classes, we need not look further. 
There cannot be more classes: all other equations yield linear combinations
of the classes. 
All of these 20 equations have real coefficients and we can use 
Lemmas 3.1-3.3 and 3.5.

Lemma 3.2 holds without any special considerations. 
In Lemma 3.3 we notice that a complex $K_3$ surface is compact and as
a real 4-manifold it is closed. Thus, the product of two $K_3$ 
surfaces is closed. $M$ is obtained from $K3\times K3$ by adding an
equation and is also closed. Lemma 3.3 can be used. 

Lemma 3.4 applies also to $M$ and gives the (2,2)-form.

The conditions for Lemma 3.5 are filled by $M$ and the lemma 
shows the (2,2)-form of Lemma 3.4 cannot 
be represented as a $\mathbb{Q}$-linear combination of classes of 
submanifolds.

The (2,2)-form in Lemma 3.4 has only one term in the local coordinates
of $M$ everywhere. Thus, it is not a linear combination. It is a single
term and therefore can be considered as a rational class.   
\end{proof}

I wrote the first version of this paper in 2011 and sent it to arXiv
for discussion purposes. I was not at all sure if my result was correct
since the solution seemed too easy and there was a conflict with a 
statement in [4] on page 54.
I got an answer from an expert of the field,  
Bert van Geemen. His comment was that the (2,2)-form in Lemma 3.4 does   
not represent a Hodge class since we cannot show that it is rational. 

The explanation by Bert van Geeman is the following:

$''$ We remember that singular cohomology of a complex
manifold of complex dimension $n$ is defined as the singular cohomology of 
the underlying real manifold of dimension $2n$. 
This yields $H^*(X;\mathbb{Z})$. In order to get $H^*(X;k)$, where the field
$k$ is $\mathbb{Q}$ or $\mathbb{R}$, or $\mathbb{C}$, 
we form the tensor product of $H^*(X;\mathbb{Z})$ with $k$.
When the base classes of $\mathbb{Z}$ are selected, the classes of the tensor
products are linear combinations of the base classes with coefficients in $k$.
We can multiply the created (2,2)-class with any number in $k$ and get a 
harmonic class. If the class is a multiple of a base vector, we can always
multiply it with a suitable number to get a class in $H^4(X;\mathbb{Q})$.
However, if it is a linear combination with generic coefficients, 
multiplication by one number does not give a rational class. This seems to
be the case with this (2,2)-form.

Theorem 2 in [4] on page 54, also published in [5], shows that the Hodge
conjecture holds for certain K3xK3 spaces, showing that in the general case 
the constructed (2,2)-form cannot be a Hodge class. 
The space of Hodge classes 
$B(H^2(S,\mathbb{Q})\otimes H^2(S,\mathbb{Q}))$ is identified up to a 
Tate twist with $End_{Hdg}(H^2(S,\mathbb{Q}))$ on page 14 in [4]. 
Zarhin's theorem is used to characterize $End_{Hdg}(T)$ on page 16, 
and Mukai's theorem
is used in the proof of Theorem 2, as in the proofs of other theorems in [4].
Mukai's theorem requires that the endomorphism $\varphi:T(S_1)\to T(S_1)$ 
is a Hodge isometry. A Hodge isometry maps $H^{2,0}(X)$ to $H^{2,0}(X)$, as is
clear e.g. from [6] page 211. 
This method in [4] seems to derive from a paper of D. Morrison [7].
The endomorphism in the presented paper is complex conjugation, which sends
$H^{2,0}(X)$ to $H^{0,2}(X)$ and is thus not a Hodge isometry. Mukai's
theorem thus cannot be used, but the identification of the space of Hodge 
classes with endomorphisms preserving the Hodge structure still holds, 
and consequently the (2,2)-form created in Lemma 3.4 does not represent
a Hodge class.$''$

I accepted this answer from an expert at that time and wrote in revised
version of the arXiv paper that the proof does not give a counterexample 
to the Hodge conjecture. However, when I have
now checked the paper, I find that the presented (2,2)-form is not a 
linear combination of (2,2)-forms. It has a single term and therefore
a multiplication with a number does give a rational form. Thus, it is 
a Hodge class. The complex manifold $M$ can be understood as a real manifold
of dimension eight. The form corresponds to a real submanifold of dimension
four that is obtained by selecting local coordinates as 
$v_1=e_{x,1}+e_{y,1}$, $v_2=e_{x,1}-e_{y,1}$, $v_3=e_{x,2}+e_{y,2}$
and $v_4=e_{x,2}-e_{y,2}$. This submanifold gives a homology class that is
the class of the submanifold. By duality there is a cohomology class 
corresponding to it. This cohomology class is inherited to singular 
cohomology of $M$ and it is not a linear combination of anything.     

The issue with the apparent conflict with [4] is avoided by noticing that
the proofs of Lemma 3.1 and 3.5 require a special condition on the K3-surface.
Therefore they do not state anything of K3xK3 spaces in the general case. 
Yet, if it follows from those results that the (2.2)-form in Lemma 3.4 
cannot be rational, then the logic of that argument should be checked.

\end{document}